\documentclass[final,twoside,11pt]{entics}
\usepackage{enticsmacro}

\usepackage{mathtools} 
\usepackage{amsmath,amssymb}
\usepackage{cleveref}

\usepackage[english]{babel}

\usepackage{tikz}
\usepackage{tikz-cd}
\usepackage{quiver}
\usepackage{booktabs} 

\usepackage{macros}
\usepackage{fixes}
\usepackage{bbm}

\def\conf{MFPS 2024}
\def\myconf{mfps40} 	
\volume{4}			
\def\volu{4}			
\def\papno{11}			
\def\mydoi{14786}%
\def\lastname{Harington, Mimram} 
\def\runtitle{Polynomials as a Kleisli category} 
\def\copynames{E. Harington, S. Mimram}    

\begin{document}
\begin{frontmatter}
  \title{Polynomials in Homotopy Type Theory\\[1ex] as a Kleisli Category}
  \author{Elies Harington\thanksref{a}}
  \author{Samuel Mimram\thanksref{b}}
  \address{LIX, CNRS, École polytechnique, Institut Polytechnique de Paris, 91120 Palaiseau, France.}
  \thanks[a]{Email: \href{samuel.mimram@polytechnique.edu}{samuel.mimram@polytechnique.edu}}
  \thanks[b]{Email: \href{harington@lix.polytechnique.fr}{harington@lix.polytechnique.fr}}
  \begin{abstract}
    Polynomials in a category have been studied as a generalization of the traditional notion in mathematics. Their construction has recently been extended to higher groupoids, as formalized in homotopy type theory, by Finster, Mimram, Lucas and Seiller, thus resulting in a cartesian closed bicategory. We refine and extend their work in multiple directions. We begin by generalizing the construction of the free symmetric monoid monad on types in order to handle arities in an arbitrary universe. Then, we extend this monad to the (wild) category of spans of types, and thus to a comonad by self-duality. Finally, we show that the resulting Kleisli category is equivalent to the traditional category of polynomials. This thus establishes polynomials as a (homotopical) model of linear logic. In fact, we explain that it is closely related to a bicategorical model of differential linear logic introduced by Melliès.

  \end{abstract}
  \begin{keyword}
    polynomial, monad, span, linear-non-linear model, linear logic, homotopy type theory
  \end{keyword}
\end{frontmatter}


\section{Introduction}
The notion of \emph{polynomial functor} in category theory, surveyed in~\cite{gambino2013polynomial}, generalizes the traditional well-known notion of polynomial in algebra. Namely, a polynomial functor is an endofunctor of sets of the form $X\mapsto\sum_{b\in B}X^{E_b}$ for some family~$(E_b)_{b\in B}$ of sets which is called the \emph{polynomial} inducing the functor. This formulation should make clear the relationship with traditional polynomials: such a functor sends an object $X$ to a (categorical) sum of monomials, which are powers of~$X$. Moreover, those have been generalized in various ways: typed variants can be considered, and one can make sense of polynomials in any locally cartesian closed category, as well as in type theory where they correspond to the notion of \emph{container}~\cite{altenkirch2015indexed}.

The category of polynomial functors over set is cartesian, but the expected closure for the cartesian product cannot be defined as easily as expected. Firstly, the right adjoint to the cartesian product does not exist for size issues: in order to have a chance to define it, one has to restrict to polynomial functors which are \emph{finitary}, \ie those which are sums of a finite number of monomials of finite degree, which were also named \emph{normal functors} by Girard~\cite{girard1988normal}.
If the category of finitary polynomial functors is indeed cartesian closed, it was observed early on by Girard and Taylor~\cite{girard1988normal,taylor1989quantitative} that it is not so in a satisfactory way: the category of finitary polynomial functors is actually a $2$-category, but the closure fails to extend as a $2$-categorical one. Moreover, the $2$-categorical structure is necessary in order to directly consider polynomials (as opposed to polynomial functors): those form a bicategory (and not a category in any reasonable way).
In order to handle this discrepancy, it was proposed by Finster et al.~\cite{finster2021cartesian}, following ideas from Kock and collaborators~\cite{kock2012data,gepner2022operads}, that this could be suitably resolved by considering finitary polynomials over groupoids instead of sets. However, the constructions cannot performed naively (for instance, it is well-known that the category of groupoids is not cartesian closed): in order for this definition to make sense, all the constructions, such as the limits involved in the definition of the composition of polynomials, have to be taken ``up to homotopy''. A definition of polynomials along these lines has been performed in the context of $2$-categories by Kock~\cite{kock2012data} and $\infty$-categories by himself, Gepner and Haugseng~\cite{gepner2022operads}, at the cost of requiring the use of quite advanced technical tools~\cite{lurie2009higher}.

It is advocated in~\cite{finster2021cartesian} that a good setting to define and study polynomials is the one of homotopy type theory~\cite{hottbook}. In this framework, types are equipped with a structure of $\infty$-groupoid~\cite{lumsdaine2010weak,van2011types} and can be interpreted as spaces~\cite{kapulkin2021simplicial}, over which one can perform constructions which are homotopy invariant by definition. In particular, the notions of finite (co)limits, which are easy to define, are actually \emph{homotopy} (co)limits. Following these ideas, polynomials have then successfully been defined in homotopy type theory and have been shown to bear a structure of cartesian closed bicategory, when restricting to finitary ones~\cite{finster2021cartesian} (up to actually truncating the category, as explained in \cref{bicategorical-model}).

In this article, we further study this model and its structure on multiple aspects. First, we show that the notion of \emph{finiteness}, which is used in order to define finitary polynomials, can actually be taken relatively to an arbitrary universe~$\V$ (which might be taken to be the one of finite types in order to recover the traditional notion). More importantly, we show here that the cartesian closed structure of the category of finitary polynomials actually comes from what we expect to be an $\infty$-categorical model of linear logic. More precisely, the category of (finitary) spans can be seen as a subcategory of the category of polynomials, which formalizes the intuition that spans are ``linear'' polynomials, and the resulting inclusion functor has a left adjoint:
\[
  \begin{tikzcd}
    \Poly\ar[r,bend left]\ar[r,phantom,"\bot"]&\ar[l,bend left]\Span
  \end{tikzcd}
\]
Moreover, this adjunction is compatible with the monoidal structures of the two categories (induced by the cartesian structure in the case of polynomials). We thus obtain what is known as a \emph{linear-non-linear   adjunction}~\cite{benton1994mixed,mellies2009categorical}, a notion which has been established as the standard categorical model of the multiplicative exponential fragment of linear logic. The adjunction induces a comonad on the category of spans, which can be understood as the exponential modality $\bang$, and we show here that the category $\Poly$ is actually the Kleisli category associated to the comonad: such a model is known as a \emph{Seely category}~\cite{seely1987linear,bierman1995categorical}.

To be more precise, the above is the situation that we expect to be holding, but showing that it is formally the case in homotopy type theory is currently out of reach because there is no known definition of $\infty$-categories in homotopy type theory~\cite{buchholtz_update_2022}. For this reason, we work here with \emph{wild categories}, which can be understood as $\infty$-categories, without requiring higher coherences. Our result should however extend to $\infty$-categories provided that we can define those. We namely expect that the constructions we provide here are actually coherent, and moreover the constructions defined by universal properties will automatically extend to the coherent setting: for instance, our result that $\Span$ is a cartesian wild category immediately implies that it is a cartesian $\infty$-category, provided that we can show that $\Span$ is an $\infty$-category.
We also explain that, by truncation, this wild category induces a (univalent) category which is a model of linear logic: this model is apparently new and of a very different nature than the traditional ones.
Moreover, we explain that the wild category can be extended into a $2$-coherent one, which induces a bicateogry, thus allowing us to recover the previous construction of~\cite{finster2021cartesian}.
This model of linear logic we obtain coincides with one constructed by Melliès in order to build a model of differential linear logic~\cite{mellies_template_2019}: the observation that the Kleisli category associated to the exponential comonad is the category of polynomials is however new. Moreover, working in homotopy type theory allows us to simplify constructions and avoid invoking arguments based on model categories~\cite{hovey2007model}.
%


\paragraph{Plan of the paper}
We begin in \cref{categories-in-hott} by recalling the general setting of homotopy type theory (\cref{homotopy-type-theory}), as well as the notion of wild category (\cref{wild-categories}) and classical structures for those (\cref{wild-structure}). In \cref{sec:span}, we define the wild category of spans and equip it with a structure of symmetric monoidal wild category with finite products and coproducts. We then construct, in \cref{exponential-comonad}, the free commutative monoid (co)monad on spans, from which we construct a structure of Seely category (\cref{span-seely}) and thus a model of linear logic (\cref{span-categorical-model}). Finally, in \cref{polynomials}, we show that the Kleisli category associated to the comonad is precisely the one of polynomials (\cref{kleisli-poly}). In \cref{related-works}, we compare with related works: the bicategory of polynomials of~\cite{finster2021cartesian} (\cref{bicategorical-model}) and template games of \cite{mellies_template_2019} (\cref{template-games}).


\section{Categories in homotopy type theory}
\label{categories-in-hott}
In this section, we recall the definition of wild categories in homotopy type theory.
We suppose the reader already familiar with the general setting and refer to the reference book~\cite{hottbook} for details.

\subsection{Homotopy type theoretic definitions and notations}
\label{homotopy-type-theory}

\paragraph{Constructions on types}
Given a type~$A$, we write $x:A$ to indicate that $x$ is an element of type~$A$.
We write $\U$ for an arbitrary universe of types, so that its elements $A : \cU$ are types.
The initial and terminal types are respectively denoted $\empty$ and $\I$. Given two types $A$ and $B$, we respectively write $A\times B$, $A\uplus B$ and $A\to B$ for their product, coproduct and function types. Given $b:B$, we write $\cst_b:A\to B$ for the constant function whose image is~$b$. Given a type $A$ and a type family \hbox{$B:A\to\U$}, the corresponding dependent sum and function types are respectively written $\dsum {x:A}{B(x)}$ and $\dprod{x:A}{B(x)}$. In the former situation, the first and second projections are respectively written $\fst:\Sigma AB\to A$ and $\snd:\dprod {x:\Sigma AB}{B(\fst(x))}$ (and similarly for the non-dependent projections $\pi_1 : A \times B \to A$ and $\pi_2 : A \times B \to B$).
%
Given types $A, B$, we write $\iota_1 : A \to A \uplus B$ and $\iota_2 : B \to A \uplus B$ for the canonical inclusions, and given a type $C$ and maps $f : A \to C$, $g : B \to C$, we write $\langle f , g \rangle$ for the induced map $A \uplus B \to C$. Given $f : A \to C, g : B \to D$, we write $f \uplus g$ for the canonical map $A \uplus B \to C \uplus D$ (\ie $f \uplus g = \langle \iota_1 \circ f , \iota_2 \circ g \rangle)$.

\paragraph{Paths}
We write $x\eqdef y$ to indicate that $x$ is defined to be~$y$, and $x\equiv y$ when $x$ is equal to $y$ by definition. The type theory also features a propositional notion of equality: given two elements $x$ and $y$ of a type~$A$, we write $x=_Ay$ (or $x=y$) for the corresponding identity type, which features a distinguished proof $\refl_x:x=x$ of reflexivity for any element~$x$. Given two identities $p:x=y$ and $q:y=z$, we write $p\pcomp q:x=z$ for the identity obtained by transitivity. The elimination principle for paths (known as \emph{path induction}) roughly states that in order to show that a property $P:\dprod{x:A}{\dprod{p:a=x}{\U}}$ holds for every possible values of~$x$ and~$p$, it is enough to show it for $x\defd a$ and $p\defd\refl_a$. By using this, it is for instance easy to show that any function $f:A\to B$ induces a function $\ap{f}:(x=y)\to(f(x)=f(y))$ defined by $\ap{f}(\refl_x)\defd \refl_{f(x)}$ witnessing that equality is a congruence.

\paragraph{Univalence}
A map $f:A\to B$ is an equivalence when it admits both a left and a right inverse. We write $A\equivto B$ for the type of equivalences between~$A$ and~$B$. An identity is always an equivalence and thus, by path induction, there is a canonical map $(A=B)\to(A\equivto B)$. The \emph{univalence axiom} states that this map is an equivalence.

\paragraph{Homotopy levels}
A type $A$ is \emph{contractible} when it satisfies $\dsum{a:A}{\dprod{x:A}{a=x}}$, which can be thought of as the fact that $A$ is a point up to homotopy. The type~$A$ is a \emph{proposition} (\resp a \emph{set}, a \emph{groupoid}) when, for every $x,y:A$, the type $x=y$ is contractible (\resp a proposition, a set). Given a type~$A$, we write $\ptrunc{A}$ for its \emph{propositional truncation}, which is the universal way of turning it into a proposition (and the set truncation $\strunc{A}$ and groupoid truncation $\gtrunc{A}$ are defined similarly), with $\trunq n-:A\to{\trunc nA}$ as canonical quotient maps. We write $\Set_\U$ (\resp $\Gpd_\U$) for the type of sets (\resp groupoids) in the universe~$\U$. Those constructions are detailed in~\cite[Chapter~7]{hottbook}.

\paragraph{Fibered / indexed equivalence}
\label{fibered-vs-indexed}
Any function $f:A\to B$ induces a type family $\fib f:B\to\U$ such that $\fib f(b)\defd\dsum{a:A}{(f(a)=b)}$ is the \emph{fiber} of~$f$ at~$b$. Conversely, any type family~$F:B\to\U$ induces a function which is the first projection $\fst:\Sigma BF\to B$ from the \emph{total space} $\Sigma BF$. These two operations form an equivalence between types over~$B$ and types families indexed by~$B$~\cite[Section~4.8]{hottbook}.

\subsection{Wild categories}
\label{wild-categories}
The notion of \emph{wild category}~\cite{capriotti_univalent_2017} can be understood as a ``non-coherent'' variant of the notion of $\infty$-category, where we do not require the expected coherences in dimension $n>1$.
In this way, wild categories can be understood as a first step toward $\infty$-categories.
They also constitute the right axiomatization when the objects of morphisms are sets~\cite[Chapter~9]{hottbook}.
As of today, coherences can also be written for $n=2$, \ie we can define bicategories, we expect that it might be possible to define tricategories following the traditional categorical definition~\cite{gordon1995coherence} but not much more, and other approaches based on complete semi-segal spaces allow writing definitions which are coherent up to an arbitrary dimension. However, whether $\infty$-categories can be defined in full generality in homotopy type theory is still an open question, and perhaps one of the most important ones~\cite{buchholtz_update_2022}. We first recall the notion of wild precategory, which is then refined into the one of wild category by further imposing a univalence property.

\begin{definition}
  A \emph{wild precategory}~$\C$ consists of type families
  \begin{align*}
    \Ob[\C]&:\cU
    &
    \Hom[\C]&: \Ob[\C] \times \Ob[\C] \to \cU
  \end{align*}
  for objects and morphisms (we omit the indices in the following), together with composition and identity operations
  \begin{align*}
    -\circ-&: \dprod{x,y,z: \Ob}{\Hom(y,z) \times \Hom(x,y) \to \Hom(x,z)}
    &
    \id&: \dprod{x: \Ob}{\Hom(x,x)}
  \end{align*}
  and coherence identities witnessing for associativity
  \begin{align*}
    \alpha&: \dprod{x,y,z,t: \Ob}{\dprod{h: \Hom(z,t)}{\dprod{g: \Hom(y,z)}{\dprod{f: \Hom(x,y)}{h \circ (g \circ f) = (h \circ g) \circ f}}}}
  \end{align*}
  as well as left and right unitality
  \begin{align*}
    \lambda&: \dprod{x,y: \Ob}{\dprod{f: \Hom(x,y)}{f \circ \id_x = f}}
    &
    \rho&: \dprod{x,y: \Ob}{\dprod{f: \Hom(x,y)}{\id_y \circ f = f}}
  \end{align*}
  For convenience, we write $x: \cC$ for $x: \Ob$ and $\cC(x,y)$ for $\Hom(x,y)$.
\end{definition}


Given a morphism $f:\C(x,y)$ for some objects $x,y:\C$, we say that~$f$ is an isomorphism when it admits both a left and a right inverse~\cite{capriotti_univalent_2017}. We write
\[
  \isIso(f)
  \defd
  \bigl(\dsum{g: \C(y,x)}{g \circ f = \id_x}\bigr)
  \times
  \bigl(\dsum{g: \C(y,x)}{f \circ g = \id_y}\bigr)
\]
for the corresponding predicate. We also write $x\isoto y$ for the type of isomorphisms in $\C(x,y)$.


\begin{lemma}
  \label[lemma]{isiso-is-prop}
  Being an isomorphism is a proposition.
\end{lemma}
\begin{proof}
  Fix a morphism $f: \cC(x,y)$ and suppose $((g,p),(h,q)): \isIso(f)$. Then we  have
  \[
    g = g \circ(f \circ h) = (g \circ f) \circ h = h
    \text.
  \]
  Hence $g$ is both a left and right inverse to $f$. Thus the type $\pa{\dsum{g: \C(y,x)}{g \circ f = \id_x}}$ is equivalent to $\dsum{g': \cC(y,x)}{g' = g}$, which is contractible, and similarly for the second component of $\isIso(f)$. We proved that $\isIso(f) \to \isContr(\isIso(f))$, which implies that $\isIso(f)$ is a proposition. Note that a similar property is shown in~\cite[Lemma 9.1.3]{hottbook} for categories, but the proof relies on the fact that hom types are sets (\cref{category}).
\end{proof}

Given a wild precategory~$\C$, the identity $\id_x$ on an arbitrary object~$x$ is always an isomorphism. We can thus define a map $\dprod{x,y:\cC}{(x = y)\to(x\isoto y)}$ by path induction such that, for $x:\C$, the image of $(x,x,\refl_x)$ is $\id_x$.

\begin{definition}
  A wild precategory is \emph{univalent} when
  the canonical map $(x = y)\to(x \isoto y)$
  is an equivalence of types for all $x, y: \cC$.
  A \emph{wild category} is a univalent wild precategory.
\end{definition}

\noindent
Note that since being an equivalence is a proposition, being univalent for a wild precategory is also a proposition.

The notion of (pre)category can be recovered as a specialization of the previous notion~\cite[Definitions~9.1.1 and 9.1.6]{hottbook}:

\begin{definition}
  \label[definition]{category}
  A \emph{precategory} (\resp \emph{category})~$\C$ is a wild precategory (\resp category) such that $\C(x,y)$ is a set for every objects $x,y:\C$.
\end{definition}

\noindent
Given a universe $\cU$, there is a wild category $\Type_\cU$ whose objects are types in $\cU$, morphisms are maps of types, composition and identities are the usual ones, and unitality and associativity hold definitionally.
%
Because of this, the above constructions generalize, to wild categories, traditional ones for types (for instance, the notion of univalence for wild precategories coincides with the usual one for types).

  %

The following construction ensures that one can always formally turn a wild category into a regular one. Note that, by univalence, the type of objects in a category has to be a groupoid~\cite[Lemma~9.1.8]{hottbook}.

\begin{definition}
  \label[definition]{wild-category-truncation}
  Given a wild precategory~$\C$, we write $\gtrunc\C$ for the precategory whose type of objects is~$\gtrunc{\Ob[\C]}$, whose type of morphisms is the map $\strunc{\Hom[\C]}:(x,y:\gtrunc{\Ob[\C]})\to\Set_\U$ induced by $\Hom[\C]$ (using the fact that $\Set_\U$ is a groupoid~\cite[Theorem~7.1.11]{hottbook}), and composition, identities and coherence laws are induced similarly by those of~$\C$ by truncation elimination.
\end{definition}

\begin{proposition}
  \label[proposition]{wild-category-truncation-univalence}
   The precategory $\gtrunc\C$ is univalent when $\C$ is.
\end{proposition}
\begin{proof}
  We write
  \[
    \pathToEquiv_{x,y}:(x =_{\Ob[\gtrunc\C]} y)\to(x \isoto y)
  \]
  for the canonical function associating to two equal objects of $\gtrunc\C$ an isomorphism between them: our aim is to show that this map is an equivalence. Since being an isomorphism is a proposition (\cref{isiso-is-prop}) and thus a groupoid~\cite[Theorem~7.1.7]{hottbook}, we can suppose given $x,y:\C$ and we have to show $(\gtrunq x=\gtrunq y)\equivto\strunc{x\isoto y}$. We know that $\gtrunq x=\gtrunq y$ is equivalent to $\strunc{x=y}$~\cite[Theorem~7.3.12]{hottbook}, which in turn is equivalent to $\strunc{x\isoto y}$ by univalence of~$\C$. Finally, we can show $\strunc{x\isoto y}\equivto(\gtrunq x\isoto\gtrunq y)$ by directly constructing an isomorphism between two types.
\end{proof}

\noindent
We expect that the above map is a left adjoint to the forgetful functor from wild (pre)categories to (pre)categories, we leave the proof of this fact for future work.

\subsubsection{Structure in wild categories}
\label{wild-structure}
Most of the usual categorical structures related to categories extend in the expected way to wild (pre)categories. In fact, many of those are developed in \cite[Chapter~9]{hottbook} in the setting of homotopy type theory for categories, but immediately generalize to wild categories. Let us list a few instances of this which will be used in the following.
\begin{itemize}
\item The cartesian product $\C\times\D$ of wild (pre)categories~$\C$ and~$\D$ is the wild (pre)category with $\Ob[\C\times\D]\defd\Ob[\C]\times\Ob[\D]$ and $(\C\times\D)((x,x'),(y,y'))\defd\C(x,y)\times\D(x',y')$ with expected composition and identities.
\item A \emph{functor} $F:\C\to\D$ between wild (pre)categories consists of a function $F:\Ob[\C]\to\Ob[\D]$ between the respective types of objects, together with functions $F_{x,y}:\C(x,y)\to\D(Fx,Fy)$ for $x,y:\C$ as well as identities $F_{x}\id_x=\id_{Fx}$ for $x:\C$ and $F(g\circ f)=Fg\circ Ff$ for composable $f$ and $g$. We can define identity and composite functors.
\item A \emph{natural transformation}~$\alpha:F\to G$ between functors $F,G:\C\to\D$ is a family $\alpha_x:\D(Fx,Gx)$ of morphisms indexed by~$x:\C$, equipped with equalities witnessing for the commutation of the expected diagrams.
\item A \emph{monad} $T:\C\to\C$ is an endofunctor equipped with natural transformations $\mu:T\circ T\To T$ and $\eta:\id_\C\To T$ together with equalities witnessing for the commutation of the expected diagrams (and comonads are defined dually).
\item A \emph{monoidal} wild precategory is a wild precategory~$\C$ equipped with an object $\I:\C$ and a functor $\otimes:\C\times\C\to\C$ together with expected natural isomorphisms (witnessing for associativity, and left and right unitality) and identities witnessing for the commutation of the two expected diagrams.
\end{itemize}
Given a comonad $(T,\delta,\varepsilon)$ on a wild precategory~$\C$, we write $\C_T$ for the Kleisli wild precategory whose objects are the same as~$\C$ and morphisms are $\C_T(x,y)=\C(Tx,y)$ with the usual notion of composition and identities~\cite[Section~VI.5]{mac2013categories}.

Structures specified by universal properties can of course also be translated to the setting of wild (pre)categories. For instance, an object $x$ of a wild precategory~$\C$ is \emph{terminal} when for every $y:\C$ the type $\C(y,x)$ is contractible. Following the usual categorical proof, two terminal objects are always isomorphic. When~$\C$ is univalent the terminal object is unique: the previous isomorphism implies that any two terminal objects are equal, \ie that the type of terminal objects is a proposition. The notion of \emph{cartesian product} of two objects $x,y:\C$ can be defined as a terminal object in the wild precategory of spans $\begin{tikzcd}[cramped,sep=small]x&\ar[l]z\ar[r]&y\end{tikzcd}$ in~$\C$. We write $\prodl:x\times y\to x$ and $\prodr:x\times y\to y$ for the two projections. Similarly, one can define the notion of \emph{pullback} of two cofinal morphisms in a wild precategory. Finally, the dual of previous notions (initial objects, coproducts, etc.) can also be defined.
Any wild category with finite products is canonically equipped with a symmetric monoidal structure, called the \emph{cartesian symmetric monoidal structure}, with the tensor product induced by cartesian product and distinguished object being the terminal object.

Our aim in this article is to show that the wild category of spans has the structure of a wild Seely category. This notion generalizes to wild categories the notion of \emph{Seely category}, in the same way as above. It was introduced in~\cite{seely1987linear} and then refined by Bierman~\cite{bierman1995categorical} who added a missing axiom, in order to provide a notion of categorical model of intuitionistic linear logic. A detailed presentation, as well as relationship with other categorical models of linear logic, can be found in~\cite{mellies2009categorical}.


\begin{definition}
  \label[definition]{seely-cat}
  A \emph{wild Seely category} is a symmetric monoidal closed wild category $(\C, \otimes, \I)$, which is moreover cartesian, with binary product and unit respectively denoted by $\with$ and $\top$, together with a comonad $(\bang, \delta, \varepsilon)$ on $\C$ called the \emph{exponential} which is symmetric monoidal, with structural natural isomorphisms $m^2_{A,B} : \bang(A \with B) \to \bang A \otimes \bang B$ and $m^0_{\top} : \bang \top \to \I$ making $\bang : (\C, \otimes, \I) \to (\C, \with, \top)$ into a symmetric monoidal functor, and such that moreover the following diagram commutes:
  \[\begin{tikzcd}[column sep=huge]
    {!(A \with B)} & {\bang \bang (A \with B)} & {\bang(\bang A \with \bang B)} \\
    {\bang A \otimes \bang B} && {\bang \bang A \otimes \bang \bang B}
    \arrow["{\delta_{A \with B}}", from=1-1, to=1-2]
    \arrow["{\delta_A\otimes\delta_B}"', from=2-1, to=2-3]
    \arrow["{\bang(\bang \pi_1, \bang \pi_2)}", from=1-2, to=1-3]
    \arrow["{m^2_{A,B}}"', from=1-1, to=2-1]
    \arrow["{m^2_{\bang A, \bang B}}", from=1-3, to=2-3]
  \end{tikzcd}\]
  The morphisms $m^2_{A,B}$ and $m^0$ are often called \emph{Seely isomorphisms}.
\end{definition}

\noindent
Any Seely category induces a linear-non-linear adjunction between itself and its Kleisli category~\cite{mellies2009categorical}.


\section{Wild categories of spans}
\label[section]{sec:span}

\begin{definition}
    Given a wild category $\cC$ with pullbacks, its wild precategory of \emph{spans}, noted $\Span(\cC)$ or $\Span_\cC$, is the wild precategory whose objects are the same as for $\cC$ and whose morphisms are spans
    $\begin{tikzcd}[cramped,sep=small]a&\ar[l,"s"']x\ar[r,"t"]&b\end{tikzcd}$,
    written $(s,t):a\spanto b$.
    The identities are
    $\begin{tikzcd}[cramped,sep=small]a&\ar[l,"\id_a"']a\ar[r,"\id_a"]&a\end{tikzcd}$
    and composition is induced by pullback:
    \begin{equation}
      \label{span-composition}
      \begin{tikzcd}[row sep=tiny]
        &&\ar[dl,dotted,"\pi_1"']x\times_b y\ar[dd,phantom,"\lrcorner"{rotate=-45,pos=0}]\ar[dr,dotted,"\pi_2"]\\
        &\ar[dl,"s"']x\ar[dr,"t"]&&\ar[dl,"u"']y\ar[dr,"v"]\\
        a&&b&&c
      \end{tikzcd}
  \end{equation}
\end{definition}

\noindent
As explained in \cref{wild-structure}, univalence entails that objects defined by universal properties are unique, so that defining composition by pullback in a category does not require any ``choice'' of pullbacks as is usually the case in set theory.
Unitality, associativity (and eventual higher coherences) of $\Span(\cC)$ also follow from the univalence of $\cC$ and the universal properties of iterated pullbacks.
Moreover, we have that
\begin{proposition}
  \label[proposition]{span-univalent}
    Given a wild category $\cC$, the wild precategory $\Span(\cC)$ is univalent.
\end{proposition}
This proposition follows easily from the gollowing lemma, together with the univalence of $\cC$.
\begin{lemma}
  \label[lemma]{iso-span-charac}
  A span $(s,t):a\spanto b$ is an isomorphism in $\Span(\cC)$ if and only if the morphisms $s$ and $t$ are isomorphisms in $\cC$.
\end{lemma}
\begin{proof}
  Suppose given a span $(s,t):a\spanto b$ which admits an inverse $(u,v):b\spanto a$. Since $(u,v)\circ(s,t)=\id_a$, with the notations of \cref{span-composition}, we have $\pi_1\circ s=\id_a$ and therefore $s$ admits a left inverse. Similarly, $s$ admits a right inverse and is thus invertible, and similarly $t$ is invertible.
\end{proof}
\begin{proof*}{Proof of \cref{span-univalent}}
  Let $a,b$ be objects of $\cC$.
  By \cref{iso-span-charac} and univalence of~$\C$, we have
  \[
    (a\equivto_{\Span}b)
    \qequivto
    \dsum{x: \cU}{(x \simeq_\cC a) \times (x \simeq_\cC b)}
    \qequivto
    \dsum{x:\cU}{(x = a)\times (x = b)}
    \qequivto
    (a=b)
  \]
  where the last step follows from classical identities and contractibility of
  singleton types~\cite[Lemma~3.11.8]{hottbook}.
\end{proof*}

\begin{remark}
  \label[remark]{rmk:selfdual}
  \label[remark]{span-self-dual}
  As evident from the definition, the category $\Span(\cC)$ is self-dual: there is an involutive contravariant functor $\Span(\cC)^\op \to \Span(\cC)$ acting as the identity on morphisms and switching the two legs of morphisms.
\end{remark}

\subsection{Functoriality of the Span construction}
The goal of this subsection is to address general question of the form: ``given some structure on $\cC$, does it lift to $\Span(\cC)$?''
An $\infty$-categorical statement that would imply all of the results we prove in this section would be along the lines of
\begin{conjecture}
  \label[conjecture]{conj:span-func}
  The $\Span$ construction underlies a limit-preserving $(\infty,2)$-functor from the $(\infty,2)$-cate\-gory of categories with pullbacks, functors preserving pullbacks, and cartesian natural transformations, to the $(\infty,2)$-category of categories.
  \Cref{extend-functor,span-nat-composition,prop:extendnattrans} together with the fact that we have the obvious equivalence \hbox{$\Span(\cC \times \cD) \equivto \Span(\cC) \times \Span(\cD)$} form a $2$-categorical approximation to this conjecture.
\end{conjecture}
\noindent
However, since we cannot state this result in full generality in homotopy type theory, let alone prove it, we detail how the span construction acts on low-dimensional cells, namely functors and natural transformations.
\begin{proposition}
  \label[proposition]{extend-functor}
  Let $\cC$ and $\cD$ be wild categories with pullbacks and $F: \cC \to \cD$ be a functor.
  If $F$ preserves pullbacks, then it lifts to a functor $\Span(F): \Span(\cC) \to \Span(\cD)$.
  $\Span(F)$ acts as $F$ on objects, and given a span $(f,g) : a \spanto b$, we have $F(f,g) \eqdef (Ff , Fg) : Fa \spanto Fb$.
\end{proposition}
\begin{proof}
  We define $\Span(F)$ on objects by
  $\Span(F)(a) \eqdef Fa$
  and on morphisms by
  \[
    \Span(F)(\begin{tikzcd}[cramped,sep=small]{a}&{}\ar[l,"f"']x\ar[r,"g"]&b\end{tikzcd})
    \defd
    \begin{tikzcd}[cramped,sep=small]Fa&\ar[l,"Ff"']Fx\ar[r,"Fg"]&Fb\end{tikzcd}
  \]
  That $\Span(F)$ preserves identities immediately follows from the fact that $F$ does, and that $F$ preserves composition of spans is precisely the statement that $F$ preserves pullbacks.
\end{proof}

Before explaining how to lift natural transformations to categories of spans, we explain how to lift morphisms.
Given objects $a, b : \C$, any morphism $f : a \to b$ canonically induces two spans
\[
  Lf : a \spanto b \eqdef \begin{tikzcd}[cramped,sep=small]a &\ar[l, "\id_a"'] a \ar[r, "f"] &b\end{tikzcd}
  \qquad\qquad \text{and} \qquad\qquad
  Rf : b \spanto a \eqdef\begin{tikzcd}[cramped,sep=small]b &\ar[l, "f"'] a \ar[r, "\id_a"] &a\end{tikzcd}
\]
Moreover, we have $L(g \circ f) = Lg \circ Lf$ since the center square in the following diagram is a pullback square:
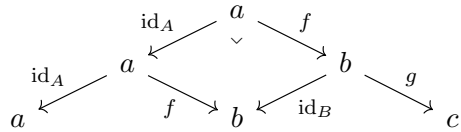
\[\begin{tikzcd}[row sep=tiny]
	&& a \\
	& a && b \\
	a && b && c
	\arrow["{\id_A}"', from=2-2, to=3-1]
	\arrow["f"', from=2-2, to=3-3]
	\arrow["{\id_B}", from=2-4, to=3-3]
	\arrow["g", from=2-4, to=3-5]
	\arrow["{\id_A}"', from=1-3, to=2-2]
	\arrow["f", from=1-3, to=2-4]
	\arrow["\lrcorner"{anchor=center, pos=0.125, rotate=-45}, draw=none, from=1-3, to=3-3]
\end{tikzcd}\]
This makes $L$ into a covariant functor $L: \cC \to \Span(\cC)$.
Similarly, $R$ extends to a contravariant functor $R: \cC^{\op} \to \Span(\cC)$. Actually, $R$ is just $L$ precomposed to the self-duality functor of \cref{rmk:selfdual}.

\begin{proposition}
    \label[proposition]{prop:extendnattrans}
    Let $F,G: \cC \to \cD$ be pullback-preserving functors between categories with pullbacks, and $\alpha: F \Rightarrow G$ be a cartesian natural transformation. Then $\alpha$ extends as a natural transformation $\Span(\alpha): \Span(F) \Rightarrow \Span(G)$.
\end{proposition}
\begin{proof}
    The natural transformation $\Span(\alpha)$ is given on objects by $\Span(\alpha)_a \eqdef L\alpha_a : Fa \spanto Ga$ for $a$ in $\Span(\cC)$.
    Now let $\begin{tikzcd}[cramped,sep=small]a&\ar[l,"f"']x\ar[r,"g"]&b\end{tikzcd}$ be a morphism in $\Span(\cC)$.
    The corresponding naturality square is the following square of spans:
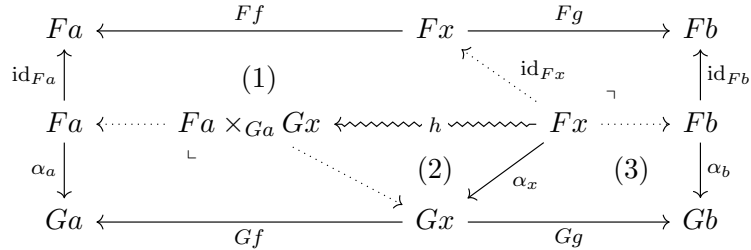
\[\begin{tikzcd}
	Fa\ar[phantom,drr,"(1)"]&& Fx && Fb \\
	Fa & {Fa \times_{Ga} Gx} &{}\ar[d,phantom,"(2)"]& \ar[dl,"\alpha_x"]Fx\ar[dr,phantom,"(3)"]& Fb \\
	Ga && Gx && Gb
	\arrow["Ff"', from=1-3, to=1-1]
	\arrow["{\id_{Fa}}", from=2-1, to=1-1]
	\arrow["{\alpha_a}"', from=2-1, to=3-1]
	\arrow["{\id_{Fb}}"', from=2-5, to=1-5]
	\arrow["{\alpha_b}", from=2-5, to=3-5]
	\arrow["Fg", from=1-3, to=1-5]
	\arrow["Gf", from=3-3, to=3-1]
	\arrow["Gg"', from=3-3, to=3-5]
	\arrow[dotted, from=2-2, to=3-3]
	\arrow[dotted, from=2-2, to=2-1]
	\arrow["{\id_{Fx}}"'{pos=0.3}, dotted, from=2-4, to=1-3]
	\arrow[dotted, from=2-4, to=2-5]
	\arrow["\lrcorner"{anchor=center, pos=0.125, rotate=-90}, draw=none, from=2-2, to=3-1]
	\arrow["\lrcorner"{anchor=center, pos=0.125, rotate=90}, draw=none, from=2-4, to=1-5]
	\arrow["h"{description}, squiggly, from=2-4, to=2-2]
      \end{tikzcd}
    \]
    where the dotted squares are pullbacks, the middle arrow~$h$ is induced by universal property of the pullback $Fa\times_{Ga}Gx$ from the naturality square
    \[
      \begin{tikzcd}
        Fa\ar[d,"\alpha_a"']&\ar[l,"Ff"']Fx\ar[d,"\alpha_x"]\\
        Ga&\ar[l,"Gf"]Gx
      \end{tikzcd}
    \]
    thus making the triangles (1) and (2) commute, and the square (3) commutes by naturality of~$\alpha$ (with respect to $g$). This can be thought of as showing that the natural transformation $\Span(\alpha):\Span(F)\to\Span(G)$ is oplax. Moreover, since $\alpha$ is supposed to be cartesian, the above naturality square is a pullback and thus $h$ is an isomorphism. We can thus conclude that $\Span(\alpha)$ is a natural transformation.
\end{proof}

\begin{proposition}
  \label[proposition]{span-nat-composition}
  The action on cartesian natural transformations $\alpha \mapsto \Span(\alpha)$ respects vertical and horizontal composition.
\end{proposition}
\begin{proof}
  First, observe that by the computation of \cref{extend-functor}, the action of $\Span(-)$ obviously preserves the composition of functors.

  First, for vertical composition of cartesian natural transformations : let $\cC, \cD$ be wild categories with pullbacks, $F,G,H : \cC \to \cD$ pullback-preserving functors,
  $\alpha : F \To G$ and $\beta : G \To H$ cartesian natural transformations.
  Given $A : \cC$, we have $\Span(\beta \circ \alpha)_A = L(\beta_A \circ \alpha_A) = L(\beta_A) \circ L(\alpha_A) = \Span(\beta_A) \circ \Span(\alpha_A)$, so $\Span$ preserves vertical composition of cartesian natural transformations.

  Now, for whiskering operations : let $\cC, \cD, \cE$ be categories with pullbacks, $F, F' : \cC \to \cD$, $G, G' : \cD \to \cE$ be pullback-preserving functors, $\alpha : F \To F'$ and $\beta : G \To G'$ be cartesian natural transformations.
  \begin{itemize}
      \item First for left whiskering : write $G \ast \alpha : G \circ F \To G \circ F'$ for the left whiskering of $\alpha$ by $G$, \ie $(G \ast \alpha)_A \eqdef G(\alpha_A)$ for all $A : \cC$.
      Since $G$ preserves pullbacks, $G \ast \alpha$ is also a cartesian natural transformation.
      Moreover, for all $A : \cC$, we have
      \[
        \Span(G \ast \alpha)_A = L(G(\alpha_A)) = G(L(\alpha_A)) = \Span(G)(\Span(\alpha)_A)
      \]
      so that $\Span$ preserves left whiskering.

      \item Then, for right whiskering : similarly, write $\beta \ast F : G \circ F \To G' \circ F$ for the right whiskering of $\beta$ by $F$, \ie $(\beta \ast F)_A = \beta_{FA}$ for all $A : \cC$.
      The naturality squares for $\beta \ast F$ are still naturality squares for $\beta$, so they are cartesian regardless of the fact that $F$ preserves pullbacks.
      Now, for all $A : \cC$, we have
      \[
        \Span(\beta \ast F)_A = L(\beta_{FA}) = \Span(\beta)_{FA} = (\Span(\beta) \ast \Span(F))_A
      \]
      so that $\Span$ preserves right whiskering.
  \end{itemize}
  We just proved $\Span$ preserves vertical composition and left and right whiskering of cartesian natural transformations, so it also preserves horizontal composition.
  This concludes the proof.
\end{proof}

Let $\cC$ be a wild category with pullbacks and a terminal object $\I$. In particular, $\cC$ admits finite products, since they are pullbacks over~$\I$. Those finite products endow $\cC$ with a symmetric monoidal structure.
\begin{proposition}
  \label[proposition]{smc-lift-span}
    The cartesian symmetric monoidal structure on $\cC$ lifts to a symmetric monoidal structure on $\Span(\cC)$.
\end{proposition}
\begin{proof}
    The cartesian product functor $-\times- : \cC \times \cC \to \cC$ commutes with pullbacks since limits commute with limits. By \cref{extend-functor}, it therefore lifts to a functor $\Span(\cC \times \cC) \to \Span(\cC)$.
    Moreover, there is a canonical equivalence of wild categories $\Span(\cC \times \cC) \simeq \Span(\cC) \times \Span(\cC)$.
    Composing $\Span(-\times-)$ with this equivalence, we obtain a tensor product $-\otimes- : \Span(\cC) \times \Span(\cC) \to \Span(\cC)$.
    The unit object of the monoidal structure is $\I$, seen as an object of $\Span(\cC)$.
    The unitality, associativity, braiding, symmetry and all higher coherences for the tensor and unit are lifted from the coherences induced by the universal properties of the cartesian product in $\cC$.
\end{proof}
\begin{remark}
    From the point of view of \cref{conj:span-func}, we are using the following idea: a symmetric monoidal $\infty$-category is a commutative monoid object in the $\infty$-category of $\infty$-categories \cite[Remark~2.4.2.6]{lurie_higher_2017}.
    We expect that finite products on $\cC$ should induce a structure of commutative monoid object on $\cC$ in the $\infty$-category of $\infty$-categories with pullbacks and pullback-preserving functors.
    By \cref{conj:span-func}, the $\infty$-functor $\Span$ should be cartesian, and thus preserves commutative monoid objects.
    Hence it should lift the cartesian symmetric monoidal structure on $\cC$ to a symmetric monoidal structure on $\Span(\cC)$.
    This idea tells us even more: any symmetric monoidal structure on $\cC$ that is compatible with pullbacks should lift to a symmetric monoidal structure on $\Span(\cC)$.
\end{remark}


\begin{proposition}
  \label[proposition]{span-monoidal-closed}
    The symmetric monoidal wild category $(\Span(\cC),\otimes,\I)$ is monoidal closed.
    For any objects $A, B$ in $\cC$, the internal hom in $\Span(\cC)$ is given by $A\llimp B \eqdef A \otimes B$.
\end{proposition}
\begin{proof}
  Given objects $A, B, C$ of $\cC$, we have
  \begin{align*}
      \Span_\cC(A \otimes B, C)
      &\equiv \Span_\cC(A \times B, C)
          &&\text{(definition of $\otimes$)}\\
      &\equivto \dsum{X \in \cC}{(X \to A \times B) \times (X \to C)}
          &&\text{(definition of $\Span_\C$)}\\
      &\equivto \dsum{X \in \cC}{(X \to A \times B \times C)}
          &&\text{(univ. prop. of $\times$ in $\cC$)}\\
      &\equivto \dsum{X \in \cC}{(X \to A) \times (X \to B \times C)}
          &&\text{(univ. prop. of $\times$ in $\cC$)}\\
      &\equivto \Span_\cC(A, B \otimes C)
          &&\text{(definition of $\Span_\C$)}
  \end{align*}
\end{proof}

\begin{remark}
  \label[remark]{span-star-autonomous}
  The wild symmetric monoidal category of spans is actually compact closed, with each object being self-dual, which can be shown following the categorical proof~\cite{stay2016compact}. As such it is monoidal closed with internal hom being given by $A\llimp B\defd A^*\otimes B\equiv A\otimes B$. Note that it also implies that it is $*$-autonomous.
\end{remark}

\subsection{Finite products in spans of types}
We now focus our attention to the case $\cC =\cU$ and show that the category $\Span_\U$ has finite products.



\begin{proposition}
  \label[proposition]{prop:span-cartesian}
  \label[proposition]{span-products}
  The wild category $\Span_\cU$ admits finite products, and they are computed as coproducts in $\cC$, \ie we have $\pi_1 = R \iota_1 : A \uplus B \spanto A$ and $\pi_2 = R \iota_2 : A \uplus B \spanto B$.
\end{proposition}
\begin{proof}
  We restrict our attention to binary products and the terminal object.
  Let $A$ and $B$ be types. We have
  \begin{align*}
      \Span_\cU(X, A \uplus B)
      &\simeq (X \times (A \uplus B)) \to \cU
          &&\text{(fibered/indexed equivalence)}\\
      &\simeq ((X \times A) \uplus (X \times B)) \to \cU
          &\\
      &\simeq ((X \times A) \to \cU) \times ((X \times B) \to \cU)
          &\\
      &\simeq \Span_\cU(X,A) \times \Span_\cU(X,B),
  \end{align*}
  which proves that $A \uplus B$ has the universal property of the cartesian product in $\Span_\cU$.

  For the terminal object, given a span $X \xleftarrow{s} Y \xrightarrow{t} \empty$, the map $t: Y \to \empty$ is necessarily unique and implies that $Y \simeq \empty$. Hence the data of such a span reduces to that of a map $\empty \to X$, of which there is only one. Hence $\empty$ is terminal in $\Span_\cU$.
\end{proof}

\begin{remark}
  The proof of \cref{span-products} could be more generally performed for any category~$\C$ which is lextensive, but our main focus will be on $\C\defd\U$ here.
\end{remark}
\begin{remark}
  \label[remark]{span-coproducts}
  By self-duality (see \cref{span-self-dual}), \cref{span-products} implies that $\Span_\U$ also has coproducts.
\end{remark}

\noindent
We now have a symmetric monoidal closed category $\Span(\cU)$, which is furthermore cartesian (\ie it admits finite products).
To enhance it to a model of linear logic following \cref{seely-cat}, we need to equip it with a comonad whose underlying functor is symmetric monoidal from the cartesian structure to the monoidal structure.


\section{The exponential comonad}
\label{exponential-comonad}
Remember that our goal is to define a comonad on $\Span(\cU)$. By self duality (\cref{rmk:selfdual}), this amounts to defining a monad on spans. We begin by defining a monad on $\Type_\cU$, and then show that it lifts to spans.

\subsection{The exponential as a monad on $\Type_\cU$}
\label{monad-on-type}
\label{small-type}
In order to axiomatize the notion of ``small type'', we suppose given a type $\V:\U$, whose elements are codes for ``small'' types, together with an embedding (\ie a fully faithful function) $\El:\V\to\U$ which to every code associates an actual type.
In the following, for simplicity, given $A:\V$, we simply write $x:A$ instead of $x:\El A$ for an element~$x$ of a small type~$A$.
Because $\El$ is an embedding, the property for $A:\U$ of being equal to $\El A'$ for some $A':\V$ is a proposition, and we say that $A$ is \emph{$\V$-small} when this is the case. We write $A\tov\U$ for the type of families $F : A \to \cU$ whose total space is $\cV$-small.
In the following, we always suppose that $\V$ is closed under dependent sums (\ie for every~$A:\V$ and $B:A\to\V$, we have $\Sigma AB:\V$), finite coproducts (in $\U$) and contains the terminal type. (Note: this axiomatization is close to the one of a \emph{regular cardinal} in set theory).

\begin{example}
  In a type theory with a hierarchy of cumulative universes $(\U_i)_{i \in \bN}$, if $\U \defd \U_j$ for some $j$, we can take $\V \defd \U_i$ for any $i < j$ to be the universe of small types.
  In that case the map $\El$ is just the inclusion of types $\U_i \hookrightarrow \U_j$.
\end{example}


\begin{example}
  \label[example]{fin-universe}
  Given a natural number $n:\N$, we write $\Fin[n]$ for a type with~$n$ elements. We define the type of \emph{finite types} as $\Fin\defd\dsum{A:\U}{\dsum{n:\N}{\ptrunc{A=\Fin[n]}}}$: an element of this type is a type together with the mere proof that it has a finite cardinal. We can then can take $\V\defd\Fin$ as universe of small types, with $\El$ being the first projection.
\end{example}


\begin{definition}
    \label[definition]{def:bang}
    Given a type $A$, its \emph{exponential} $\bang_\cV A$ (read ``bang A'') relative to $\cV$ is defined to be
    \[\bang_\cV A \eqdef \dsum{E: \cV}{(E \to A)}\]
    We often simply write $\bang A$, omitting $\cV$ when it is obvious from the context.
    The exponential acts on morphisms by postcomposition: given a map $f : A \to B$, the induced map $\bang f : \bang A \to \bang B$ is defined by $\bang f(E,p) \eqdef (E, f \circ p)$.
    This makes $\bang$ into an endofunctor on $\Type_\cU$.
\end{definition}

\begin{remark}
    \label[remark]{rmk:pov-change}
    Consider the case $\cV \equiv \Fin$ of \cref{fin-universe}.
    Then $\bang A$ is the homotopical analogue of the multiset construction.
    In ordinary set theory, multisets on a set $A$ are usually defined as maps $A \to \bN$ with finite support.
    In other words, a multiset on $A$ is a finite collection of elements of $A$, with possible repetitions.
    Another way to think of a multiset $X$ on $A$ is simply as a finite set $X$ whose elements are \emph{colored} by elements of $A$.
    This intuition is implemented formally by alternatively defining a multiset on $A$ to be the data of a finite set $X$ together with a map $c: X \to A$.
    It is this latter intuition that we chose as \ref{def:bang} here.

    Under the fibered/indexed equivalence (\cref{fibered-vs-indexed}), we could have also adapted the ``map with finite support'' point of view, defining a multiset on $A$ to be a map $F : A \tov \U$, \ie a map $F:A\to\U$ such that its total space $\dsum{x : A}{F(x)}$ is a $\cV$-small type (notice how ``finite support'' has to be replaced with ``finite total space'' to make sense for an arbitrary base type $A$).

    We say that $\bang A$ is a \emph{homotopical} version of multisets in the sense that, contrary to its set-theoretic version, it remembers more symmetries.
    For instance, consider the set $A = \{a,b\}$, and the multiset $X = [a,a,b]$.
    Seen as an element of $\bang A$, $X$ has two self-identifications: the identity, and the one that swaps the two elements \emph{colored} by $a$.
    Hence the type $X =_{\bang A} X$ has two elements, and $\bang A$ is a groupoid that is not a set, even though $A$ was.
\end{remark}

\begin{remark}
    Consider still $\cV \equiv \Fin$.
    Given a set $A$, the set of multisets on $A$ is the free commutative monoid on $A$.
    Similarly, given a type $A$ seen as an $\infty$-groupoid, we expect that $\bang A$ is the free symmetric monoidal $\infty$\nbd-groupoid on $A$ (and similarly we expect $\bang_\cV A$ to be the free symmetric monoid with $\cV$-small sums).
\end{remark}

In order to explore the monad laws of $\bang$, we can first informally use the analogy with multisets. Writing $\Mul(A)$ for the multiset on a set~$A$, we expect the unit $\eta_A:A\to\Mul(A)$ to be the map sending $a$ to the multiset with $a$ (once) as only element. Similarly, we expect the multiplication $\mu_A:\Mul(\Mul(A))\to\Mul(A)$ to take a multiset of multisets to the multiset of its elements.
%
This suggests defining the following natural transformations
\begin{align*}
  \eta_A:A & \to\bang A & \mu_A:\bang \bang A & \to \bang A\\
  a & \mapsto (\I, \cst_a)&(E, p) & \mapsto \Bigl(\bigl(\dsum{e: E}{\fst(p(e))}\bigr), \pa{(e,e') \mapsto \snd(p(e))(e')}\Bigr)
\end{align*}
Note that the above definitions use the fact that $\V$ contains the terminal object and is closed under dependent sums.

\begin{remark}
    \label[remark]{rmk:simplifymu}
    By definition, we have $\bang \bang A \equiv \dsum{E: \cV}{\pa{ E \to \dsum{E': \cV}{\pa{E' \to A}}}}$, so an element in $\bang \bang A$ consists of a $\cV$-small type $E$, and a family $E'$ of $\cV$-small types indexed by $E$, together with a map associating to each pair $(e,e'): \sum_{e:E} E'(e)$ a ``color'' in~$A$.
    In other words, $\bang \bang A \equivto \dsum{E : \cV}{\dsum{E' : E \to \cV}{(\Sigma E E' \to A)}}$.
    Under this equivalence, the multiplication map $\mu: \bang \bang A \to \bang A$ simply sends the triplet $(E,E',p)$ to the total space $\Sigma E E'$ colored by $p$.
    Throughout the rest of this text, we will often abuse notation by writing $(E,E',p) : \bang \bang A$, implicitly using this equivalence.
\end{remark}

\begin{lemma}
  \label[lemma]{eta-mu-natural}
  The transformations $\eta$ and $\mu$ defined above are natural.
\end{lemma}
\begin{proof}
  Let $f : A \to B$ be a map of types and consider the follow naturality squares:
  \[\begin{tikzcd}
      A & {\bang A} && {\bang \bang A} & {\bang A} \\
      B & {\bang B} && {\bang \bang B} & {\bang B}
      \arrow["{\eta_B}"', from=2-1, to=2-2]
      \arrow["{\eta_A}", from=1-1, to=1-2]
      \arrow["f"', from=1-1, to=2-1]
      \arrow["{\bang f}", from=1-2, to=2-2]
      \arrow["{\mu_A}", from=1-4, to=1-5]
      \arrow["{\bang \bang f}"', from=1-4, to=2-4]
      \arrow["{\mu_B}"', from=2-4, to=2-5]
      \arrow["{\bang f}", from=1-5, to=2-5]
  \end{tikzcd}\]
Given $a : A$, we have
\[
  (\bang f)(\eta_A(a)) = (\bang f)(\I, \cst_a) = (\I, \cst_{f(a)}) = \eta_B(f(a))
\]
(and this equality is actually reflexivity).

Given $(E,E',p) : \bang \bang A$ (under the equivalence of \cref{rmk:simplifymu}), we have
\[
  (\bang f)(\mu_A(E,E',p)) = (\bang f)(\Sigma E E', p) = (\Sigma E E', f \circ p)
\]
and
\[
  \mu_B((\bang \bang f)(E,E',p)) = \mu_B(E,E',f\circ p) = (\Sigma E E', f \circ p)
\]
which completes the proof of naturality.
\end{proof}

\begin{proposition}
    \label[proposition]{bang-monad}
    The triple $(\bang, \mu, \eta)$ is a monad on $\U$.
\end{proposition}
\begin{proof}
  Let $A : \U$. We need to prove the following square and triangles commute
  \[\begin{tikzcd}
      {\bang \bang \bang A} & {\bang \bang A} & {\bang A} & {\bang \bang A} & {\bang A} \\
      {\bang \bang A} & {\bang A} && {\bang A}
      \arrow["{\bang\mu_A}", from=1-1, to=1-2]
      \arrow["{\mu_{\bang A}}"', from=1-1, to=2-1]
      \arrow["{\mu_A}"', from=2-1, to=2-2]
      \arrow["{\mu_A}", from=1-2, to=2-2]
      \arrow["{\bang\eta_A}", from=1-3, to=1-4]
      \arrow["{\eta_{\bang A}}"', from=1-5, to=1-4]
      \arrow["\mu"', from=1-4, to=2-4]
      \arrow["{\id_{\bang A}}"', from=1-3, to=2-4]
      \arrow["{\id_{\bang A}}", from=1-5, to=2-4]
  \end{tikzcd}\]
  Let $(E,p) : \bang A$.
  We have
  \[
    \mu_A(\eta_{\bang A}(E,p)) \equivto \mu_A(\I,\cst_E,(*,e) \mapsto p(e)) \equivto (E, p)
  \]
  and
  \[
    \mu_A(\bang \eta_A(E,p)) \equivto \mu_A(E, e \mapsto \I, (e,*) \mapsto p(e)) \equivto (E,p)
  \]
  hence the triangles commute.

  Now let $(E,E',E'',p) : \bang \bang A$. Here we are using a characterization of $\bang \bang \bang A$ similar to that of $\bang \bang A$ in \cref{rmk:simplifymu}.
  In other words, we have $E : \V$, $E' : E \to \V$, $E'' : \Sigma E E' \to \V$ and $p : \Sigma E (\Sigma E' E'') \to A$.
  We have
  \[
    \mu_A(\bang\mu_A(E,E',E'',p)) = \mu_A(E, \Sigma E' E'', p) = (\Sigma E (\Sigma E' E'') , p)
  \]
  and
  \[
    \mu_A(\mu_{\bang A}(E,E',E'',p)) = \mu_A(\Sigma E E', E'', p) = \mu_A(\Sigma E (\Sigma E' E''), p)
  \]
  hence the square commutes, which completes the proof.
\end{proof}

\noindent
We finally construct two isomorphisms in $\cU$ that will become the Seely isomorphisms for $\Span(\cU)$ in the next section.
\begin{proposition}
    \label[proposition]{seely-isos}
    We have natural isomorphisms $l^2_{A,B}: \bang A \times \bang B \xrightarrow{\sim} \bang(A \uplus B)$ for all types $A, B: \cU$, and an isomorphism $l^0 : \I \xrightarrow{\sim} \bang \empty$. The former is defined by
    $
    l^2_{A,B}((E,p),(F,q))\defd (E \uplus F, p \uplus q)
    $.
\end{proposition}
\begin{proof}
  For $l^0$, we have
  $
  \bang \empty \equivto(\empty \tov \cU) \equivto \I
  $
  by the universal property of the empty type and the fact that it is $\V$-small.
  For $l^2$, given types $A, B: \cU$, we have
  \[
    \bang A \times \bang B
    \equivto (A \tov \cU) \times (B \tov \cU)
    \equivto (A \uplus B \tov \cU)
    \equivto \bang(A \uplus B)
  \]
  using \cref{rmk:pov-change} and closure of~$\V$ under coproducts.
  Unfolding this equivalence, we have
  \[
    l^2_{A,B}((E,p)(F,q)) \eqdef (E \uplus F, p \uplus q)
    \text.
  \]
  Given maps $f : A \to C$ and $g : B \to D$, we need to show the following naturality square commutes. 
    \[\begin{tikzcd}
        {!A \times !B} & {!(A \sqcup B)} \\
        {!C \times !D} & {!(C \sqcup D)}
        \arrow["{!(f \sqcup g)}", from=1-2, to=2-2]
        \arrow["{!f \times !g}"', from=1-1, to=2-1]
        \arrow["{l^2_{A,B}}", from=1-1, to=1-2]
        \arrow["{l^2_{C,D}}"', from=2-1, to=2-2]
    \end{tikzcd}\]
    Let $((E,p),(F,q)) : \bang A \times \bang B$.
    We have
    \begin{align*}
        (\bang(f \uplus g))(l^2_{A,B}((E,p),(F,q))) &=
            (\bang(f \uplus g))(E \uplus F, p \uplus q) =
            (E \uplus F, f \circ p \uplus g \circ q)\\
        l^2_{C,D}((\bang f \times \bang g)((E,p),(F,q))) &=
            l^2_{C,D}((E, f \circ p), (F, g \circ q)) = 
            (E \uplus F,  f \circ p \uplus g \circ q)
    \end{align*}
    which concludes the proof.
\end{proof}


\subsection{Lifting the exponential monad to spans}
We now show that the monad $\bang$ on~$\U$ lifts to a monad on~$\Span(\U)$.


\begin{proposition}
  \label[proposition]{bang-lift-span}
  The functor $\bang: \cU \to \cU$ preserves pullbacks and therefore, by \cref{extend-functor}, lifts to a functor $\Span(\bang) : \Span_\cU \to \Span_\cU$.
\end{proposition}
\begin{proof}
  Using \cref{extend-functor}, we only need to show that $\bang : \cU \to \cU$ preserves pullbacks.
  Let $A \xrightarrow{f} C \xleftarrow{g} B$ be a diagram of types.
  We have the following chain of equivalences:
  \allowdisplaybreaks
  \begin{align*}
      \bang A \times_{\bang C} \bang B
      &\equiv \sum_{E: \cV}\sum_{p: E \to A}\sum_{F: \cV}\sum_{q: F \to B} ((E, f \circ p) = (F, g \circ q))&&
          \text{(def. of $\bang$ and pullback)}\\
      &\simeq \sum_{E: \cV}\sum_{p: E \to A}\sum_{F: \cV}\sum_{q: F \to B} \Bigpa{\sum_{l: E \simeq F} (f \circ p = g \circ q \circ l)}&&
          \text{(equality in $\Sigma$-types and univalence)}\\
      &\simeq \sum_{E: \cV}\sum_{F: \cV}\sum_{l: E \simeq F}\sum_{p: E \to A}\sum_{q: F \to B} (f \circ p = g \circ q \circ l)&&
          \text{(reordering of terms)}\\
      &\simeq \sum_{E: \cV}\sum_{p: E \to A}\sum_{q: E \to B} (f \circ p = g \circ q)&&
          \text{(contracting $(F,l)$ onto $(E, \id_E)$)}\\
      &\simeq \sum_{E: \cV}\sum_{p: E \to A}\sum_{q: E \to B} \Bigpa{\prod_{e:E} (f(p(e)) = g(q(e)))}&&
          \text{(function extensionality)}\\
      &\simeq \sum_{E: \cV} \pa{E \to A \times_C B }&&
          \text{(univ. prop. of the pullback)}\\
      &\equiv \bang(A \times_C B)&&
          \text{(def. of $\bang$)}
  \end{align*}
\end{proof}

\noindent
In the following, to make notations more readable, we will still denote by $\bang$ its lifting $\Span(\bang)$.

Our aim is now to show that this functor inherits monad laws from the ones of the monad constructed on spans in \cref{monad-on-type}.
Using propositions \ref{prop:extendnattrans} and \ref{span-nat-composition}, we only need to show that the natural transformations $\eta$ and $\mu$ are cartesian.

\begin{proposition}
  \label[proposition]{eta-cartesian}
  The natural transformation $\eta: \id_{\cU} \Rightarrow \bang$ is cartesian.
\end{proposition}
\begin{proof}
  Let $f: A \to B$ be a map between types.
  Unwinding the corresponding naturality square of $\eta$, we get
\[\begin{tikzcd}
A &&& B \\
\\
{\sum_{E:\cV}(E \to A)} &&& {\sum_{E:\cV}(E \to B)}
\arrow["f", from=1-1, to=1-4]
\arrow["{a \mapsto(\top,* \mapsto a)}"{description}, from=1-1, to=3-1]
\arrow["{(E,p) \mapsto (E, f \circ p)}"', from=3-1, to=3-4]
\arrow["{b \mapsto(\top,* \mapsto b)}"{description}, from=1-4, to=3-4]
\end{tikzcd}\]
  That square is definitionally commutative, i.e.\ the left-bottom and top-right compositions compute to the same term, or still in other terms, the witness of commutativity is $\refl_{\pa{a \mapsto \pa{\top,\cst_{f(a)}}}}$.

  We have the following chain of equivalences:
  \allowdisplaybreaks
  \begin{align*}
      \bang A \times_{\bang B} B
      &\equiv \sum_{E: \cV} \sum_{p: E \to A} \sum_{b: B} \pa{(E,f \circ p) = (\top , \cst_b)}
          &&\text{(unfolding the definitions)}\\
      &\simeq \sum_{E: \cV} \sum_{l: E \simeq \top} \sum_{p: E \to A} \sum_{b:B} (f \circ p = \cst_b \circ l)
          &&\text{(equality in $\Sigma$-types, univalence, reordering)}\\
      &\simeq \sum_{p: \top \to A} \sum_{b: B} (f \circ p = \cst_b)
          &&\text{(contracting $(E,l)$ onto $(\top, \id_\top)$)}\\
      &\simeq \sum_{a: A} \sum_{b: B} (f(a) = b)
          &&\text{(universal prop. of $\top$)}\\
      &\simeq A
          &&\text{(contracting the last two terms onto $(f(a),\refl)$)}
  \end{align*}
  Moreover, computing this chain of equivalences from bottom to top exactly gives the connecting map $A \to \bang A \times_{\bang B} B$ of the pullback of the naturality square of $\eta$.
  So $\eta$ is indeed cartesian.
\end{proof}

\begin{proposition}
  \label[proposition]{mu-cartesian}
  The natural transformation $\mu: \bang\bang \Rightarrow \bang$ is cartesian.
\end{proposition}
\begin{proof}
  \begin{align*}
      {\blue \bang A} {\red \times_{\bang B}} {\green \bang \bang B}
      &\equivto
          {\blue \dsum{E : \V} \dsum{p : E \to A}}
          {\green \dsum{F : \V} \dsum{F' : F \to \V} \dsum{q : \Sigma F F' \to B}}
          {\red (E , f \circ p) =_{\bang B} (\Sigma F F' , q)}
          &&\text{(see \cref{rmk:simplifymu} for $\green \bang \bang B$)}\\
      &\equivto
          {\blue \dsum{E : \V} \dsum{p : E \to A}}
          {\green \dsum{F : \V} \dsum{F' : F \to \V} \dsum{q : \Sigma F F' \to B}}
          {\red \dsum{l : E \equivto \Sigma F F'} f \circ p = q \circ l}
          &&\text{(equality in $\Sigma$-types and univalence)}\\
      &\equivto
          {\green \dsum{F : \V} \dsum{F' : F \to \V}}
          {\blue \dsum{p : \Sigma F F' \to A}} {\green \dsum{q : \Sigma F F' \to B}}
          {\red f \circ p = q}
          &&\text{(contracting $({\blue E}, {\red l})$ onto $(\Sigma F F', \id)$)}\\
      &\equivto
          {\green \dsum{F: \cV}\dsum{F': F \to \cV}}{\blue \pa{\Sigma F F' \to A }}
          &&\text{(contracting $({\green q}, {\red -})$ onto $(f \circ p, \refl)$)}\\
      &\equivto \bang \bang A
          &&\text{(\cref{rmk:simplifymu})}
  \end{align*}
  Like with $\eta$, one can check that the underlying map of the inverse of the equivalence thus constructed is equal to the map $\bang \bang A \to \bang A \times_{\bang B} \bang \bang B$ induced by the naturality square of $\mu$.
  Hence $\mu$ is cartesian.
\end{proof}

Since $\eta$ and $\mu$ are cartesian, they lift to natural transformations in spans, making $\bang$ into a monad on $\Span(\cU)$.
By \cref{rmk:selfdual}, they also make bang into a comonad on $\Span(\cU)$.
We write $\epsilon$ and $\delta$ respectively for the counit and comultiplication of the comonad $\bang$.
Unfolding the definition, we have
\begin{align*}
    \epsilon_A &: \bang A \spanto A & \delta_A &: \bang A \spanto \bang \bang A\\
    \epsilon_A & \eqdef R(\eta_A) = \bang A \xleftarrow{\eta_A} A \xrightarrow{\id_A} A
    &
    \delta_A & \eqdef R(\mu_A) = \bang A \xleftarrow{\mu_A} \bang \bang A \xrightarrow{\id_{\bang \bang A}} \bang \bang A
\end{align*}
Similarly, $l^2$ being a natural isomorphism, its naturality squares are cartesian so that, under the equivalence $\Span(\U\times\U) \equivto \Span(\U) \times \Span(\U)$ and by selfduality (\cref{rmk:selfdual}), it lifts to a natural isomorphism $m^2_{A,B} \eqdef R(l^2_{A,B}) : \bang(A \uplus B) \spanto \bang A \times \bang B$.
The morphism $l^0$ also lifts to a morphism $m^0 \eqdef R(l^0) : \bang \empty \spanto \I$.

\begin{theorem}
  \label[theorem]{span-seely}
  The symmetric monoidal closed, cartesian wild category $(\Span(\U), \otimes, \I, \uplus, \empty)$, equipped with the comonad $(\bang,\delta,\epsilon)$ and the morphisms $m^2, m^0$, is a wild Seely category.
\end{theorem}
\begin{proof}
    We need to check that $(\bang, m^2, m^0) : (\Span(\U), \uplus, \empty) \to (\Span(\U), \otimes, \I)$ is symmetric monoidal, and that for all $A, B : \U$ the following diagram commutes
    \[\begin{tikzcd}[column sep=large]
        {\bang (A \uplus B)} & {\bang \bang (A \uplus B)} & {\bang(\bang A \uplus \bang B)} \\
        {\bang A \times \bang B} && {\bang \bang A \times \bang \bang B}
        \arrow["{\delta_A \times \delta_B}"', "\shortmid"{marking}, from=2-1, to=2-3]
        \arrow["{\delta_{A \uplus B}}", "\shortmid"{marking}, from=1-1, to=1-2]
        \arrow["{\bang(\bang \pi_1, \bang \pi_2)}", "\shortmid"{marking}, from=1-2, to=1-3]
        \arrow["{m^2_{A,B}}"', from=1-1, to=2-1]
        \arrow["{m^2_{\bang A, \bang B}}", from=1-3, to=2-3]
    \end{tikzcd}\]
    But this diagram is precisely the image by the functor $R : \U^\op \to \Span(\U)$ of the diagram
    \[\begin{tikzcd}[column sep=large]
        {\bang (A \uplus B)} & {\bang \bang (A \uplus B)} & {\bang(\bang A \uplus \bang B)} \\
        {\bang A \times \bang B} && {\bang \bang A \times \bang \bang B}
        \arrow["{\mu_A \times \mu_B}", from=2-3, to=2-1]
        \arrow["{\bang(\bang \iota_1, \bang \iota_2)}"', from=1-3, to=1-2]
        \arrow["{\mu_{A \uplus B}}"', from=1-2, to=1-1]
        \arrow["{l^2_{A,B}}", from=2-1, to=1-1]
        \arrow["{l^2_{\bang A, \bang B}}"', from=2-3, to=1-3]
    \end{tikzcd}\]
    Let $\bE \eqdef (E,E',p) : \bang \bang A$ and $\bF \eqdef (F,F',q) : \bang \bang B$ (under the equivalence of \cref{rmk:simplifymu}).
    We have
    \[l^2_{A,B}((\mu_A \times \mu_B)(\bE, \bF)) =
        l^2_{A,B}((\Sigma E E', p),(\Sigma F F', q)) =
        (\Sigma E E' \uplus \Sigma F F',  p \uplus q )\]
    and
    \[\begin{tikzcd}[sep=small]
        {\bang \bang A \times \bang \bang B} & {\bang(\bang A \uplus \bang B)} & {\bang \bang (A \uplus B)} & {\bang (A \uplus B)} \\
        {(\bE , \bF)} & {(E \uplus F, (E',p) \uplus (F',q))} & {(E \uplus F,\langle E', F' \rangle ,  p \uplus q)} & {(\dsum{x : E \uplus F} \dsum{y : \langle E', F' \rangle(x)} (p \uplus q)(x,y))}
        \arrow["{\bang\langle\bang \iota_1, \bang \iota_2\rangle}", from=1-2, to=1-3]
        \arrow["{\mu_{A \uplus B}}", from=1-3, to=1-4]
        \arrow["{l^2_{\bang A, \bang B}}", from=1-1, to=1-2]
        \arrow[maps to, from=2-1, to=2-2]
        \arrow[maps to, from=2-2, to=2-3]
        \arrow[maps to, from=2-3, to=2-4]
    \end{tikzcd}\]
    Those two elements of $\bang(A \uplus B)$ are equal by virtue of $\Sigma$-types distributing over disjoint sums.

    \medskip
    \noindent
    What remains to be shown is that the following data constitutes a symmetric monoidal functor:
  \[
    (\bang, m^2, m^0) : (\Span(\U), \uplus, \empty) \to (\Span(\U), \times, \I)
    \text.
  \]
  Since the functor $\bang$ on spans is obtained as a lifting of $\bang : \U \to \U$, and similarly for the natural transformations $m^2$ and $m^0$, using \cref{span-nat-composition} we only need to show that
  \[(\bang, l^2, l^0) : (\U, \uplus, \empty) \to (\U, \times, \I)\]
  is a symmetric monoidal functor.
  Following \cite[Section 7.3]{mellies2009categorical} (and writing respectively $\alpha, \lambda, \rho, \gamma$ for the obvious associativity, left-unitality, right-unitality, and symmetry isomorphisms), this amounts to showing the following four diagrams commute.
  \[\begin{tikzcd}
      {(\bang A \times \bang B) \times \bang C} & {\bang A \times (\bang B \times \bang C)} && {\bang A \times \I} & {\bang A} \\
      {\bang(A \uplus B) \times \bang C} & {\bang A \times \bang(B \uplus C)} && {\bang A \times \bang \empty} & {\bang(A \uplus \empty)} \\
      {\bang ((A \uplus B) \uplus C)} & {\bang(A \uplus (B \uplus C))} && {\I \times \bang B} & {\bang B} \\
      {\bang A \times \bang B} & {\bang B \times \bang A} && {\bang \empty \times \bang B} & {\bang(\empty \uplus B)} \\
      {\bang(A \uplus B)} & {\bang(B \uplus A)}
      \arrow["{\alpha_\times}", from=1-1, to=1-2]
      \arrow["{l^2_{A,B}\times\id_{\bang C}}"', from=1-1, to=2-1]
      \arrow["{l^2_{A \uplus B, C}}"', from=2-1, to=3-1]
      \arrow["{\bang\alpha_\uplus}"', from=3-1, to=3-2]
      \arrow["{\id_{\bang A} \times l^2_{B,C}}", from=1-2, to=2-2]
      \arrow["{l^2_{A, B \uplus C}}", from=2-2, to=3-2]
      \arrow["{(1)}"{description}, draw=none, from=2-1, to=2-2]
      \arrow["{\rho_\times}", from=1-4, to=1-5]
      \arrow[""{name=0, anchor=center, inner sep=0}, "{\id_{\bang A} \times l^0}"', from=1-4, to=2-4]
      \arrow["{l^2_{A,\empty}}"', from=2-4, to=2-5]
      \arrow[""{name=1, anchor=center, inner sep=0}, "{\bang\rho_\uplus}"', from=2-5, to=1-5]
      \arrow["{\lambda_\times}", from=3-4, to=3-5]
      \arrow[from=4-4, to=4-5]
      \arrow[""{name=2, anchor=center, inner sep=0}, "{\bang\lambda_\uplus}"', from=4-5, to=3-5]
      \arrow[""{name=3, anchor=center, inner sep=0}, "{l^0 \times \id_{\bang B}}"', from=3-4, to=4-4]
      \arrow["{\gamma_\times}", from=4-1, to=4-2]
      \arrow[""{name=4, anchor=center, inner sep=0}, "{l^2_{B,A}}", from=4-2, to=5-2]
      \arrow[""{name=5, anchor=center, inner sep=0}, "{l^2_{A,B}}"', from=4-1, to=5-1]
      \arrow["{\bang\gamma_\uplus}"', from=5-1, to=5-2]
      \arrow["{(4)}"{description}, draw=none, from=5, to=4]
      \arrow["{(2)}"{description}, draw=none, from=0, to=1]
      \arrow["{(3)}"{description}, draw=none, from=3, to=2]
  \end{tikzcd}\]
  Most verifications are straightforward:
  \begin{itemize}
      \item for diagram $(1)$, given $(((E,p),(F,q)),(G,r)) : (\bang A \times \bang B) \times \bang C$, both paths in the hexagon evaluate to $(E \uplus (F \uplus G), p \uplus (q \uplus r)) : \bang(A \uplus (B \uplus C))$ up to associativity of $\uplus$, which is an equivalence and hence an equality by univalence (see the case of $(4)$ for more a detailed reasoning along those lines),
      \item diagrams $(2)$ and $(3)$ are trivial,
      \item and finally for diagram $(4)$, given $((E,p),(F,q)) : \bang A \times \bang B$, we have
          \[l^2_{B,A}(\gamma_\times((E,p),(F,q))) = l^2_{B,A}((F,q),(E,p)) = (F \uplus E, q \uplus p)\]
          and
          \[\bang\gamma_\uplus(l^2_{A,B}((E,p),(F,q))) = \bang\gamma_\uplus(E \uplus F, p \uplus q) = (E \uplus F, \gamma_\uplus \circ (p \uplus q))\]
          By characterization of equality in $\Sigma$-types and univalence, an equality between $(E \uplus F, \gamma_\uplus \circ (p \uplus q))$ and $(F \uplus E, q \uplus p)$ in $\bang(B \uplus A)$ consists of the data of an equivalence $f : E \uplus F \xrightarrow{\equivto} F \uplus E$ and an equality $e : \gamma_\uplus \circ (p \uplus q) = (q \uplus p)\circ f$ in $E \uplus F \to B \uplus A$.
          Choosing $f \eqdef \gamma_\uplus$, a simple computation shows the following square commutes
          \[\begin{tikzcd}
              {E \uplus F} & {A \uplus B} \\
              {F \uplus E} & {B \uplus A}
              \arrow["{p \uplus q}", from=1-1, to=1-2]
              \arrow["{\gamma_\uplus}"', from=1-1, to=2-1]
              \arrow["{q \uplus p}"', from=2-1, to=2-2]
              \arrow["{\gamma_\uplus}", from=1-2, to=2-2]
          \end{tikzcd}\]
          which concludes the proof.
  \end{itemize}
\end{proof}

\noindent
The above theorem can be interpreted as the fact that $\Span(\U)$ is a ``wild model'' of intuitionistic linear logic. Since it is $*$-autonomous (\cref{span-star-autonomous}), it also extends to a model of classical linear logic and the presence of products (\cref{span-products}) and coproducts (\cref{span-coproducts}) allow for modeling additives.

It is not difficult to show that the categorical structures are preserved by truncation (\cref{wild-category-truncation}), and preservation of univalence was shown in \cref{wild-category-truncation-univalence}. We thus have:

\begin{theorem}
  \label[theorem]{span-categorical-model}
  The category $\gtrunc{\Span(\U)}$ is a Seely category and thus a model of linear logic.
\end{theorem}

\noindent
Note that the model we obtain, which is apparently new, is of a very different nature than the usual ones for linear logic (relational model, coherence spaces, etc.). Namely, the category $\gtrunc{\Type_\U}$ of types and homotopy classes of functions embeds, via the functor induced by~$L$, into $\gtrunc{\Span(\U)}$. This former category can be interpreted as the homotopy category of spaces, and it is known that it is not \emph{concrete}~\cite{freyd1970homotopy}, \ie cannot be recovered as a subcategory of~$\Set$, thus neither can be our model. We leave for future work the investigation of possible more direct descriptions of this model.


\section{The exponential modality and polynomials}
\label{polynomials}

We show here that the wild Kleisli category associated to our comonad on $\Span(\U)$ is the well-known category of polynomials (up to homotopy).

\begin{definition}
  \label[definition]{polynomial}
    Let $I, J$ be types. A polynomial from $I$ to $J$ is a diagram of the shape
    \begin{equation}
      \label{polynomial-diagram}
      \begin{tikzcd}
        I&\ar[l,"s"']E\ar[r,"p"]&B\ar[r,"t"]&J
      \end{tikzcd}
    \end{equation}
    %
    In other words, the type of polynomials with source $I$ and target $J$ is
    \[
      \Poly(I,J) \eqdef \dsum{E: \cU}{\dsum{B: \cU}{(E \to I) \times (E \to B) \times (B \to J)}}
    \]
\end{definition}

\noindent  
To make sense of this definition, one must understand what the data of a polynomial represents, namely, a \emph{polynomial functor}.

\subsection{The functor associated to a polynomial}
In the following, given a map $f:A\to B$ and $b:B$, we write $A_b$ for $\fib f(b)$, the function being often implicit from the context. A polynomial as in \cref{polynomial} can be thought of as the data describing a (colored) polynomial functor in the following sense. Given $j:J$, $B_j$ is the type of monomials colored by~$j$; given $b:B_j$, the arity of the monomial~$b$ is $E_b$; finally, the function $s : E_b\to I$ associates, to each variable $e:E_b$ a color in $I$.
This interpretation suggests associating to any polynomial a function, following the now classical definition of a polynomial functor~\cite{gambino2013polynomial}.

\begin{definition}
    \label[definition]{def:poly-functor}
    Any polynomial $P = (E,B,s,t,p)$ as in \cref{polynomial-diagram} induces a map $F_P: \cU^I \to \cU^J$, defined by
    \[
      F(X)(j) \eqdef \dsum{b:B_j}{\dprod{e \in E_b}{X(s(e))}}
    \]
    called the \emph{polynomial functor} induced by $P$.
\end{definition}

\noindent
Above, we abusively use the term ``functor'' following the traditional terminology which comes from the fact that, in category theory, a polynomial in a locally cartesian closed category $\cC$ between objects $x$ and $y$ induces a functor $\cC_{/x} \to \cC_{/y}$. Here, we have $\cC \eqdef \U$, and under the fibered/indexed equivalence the slice category $\U_{/X}$ is equivalent to $\U^X$.

\begin{definition}
  \label[definition]{def:linear-poly}
  A polynomial $P = (E,B,s,p,t)$ is said to be \emph{linear} if the map $p$ is an equivalence.
  This is precisely equivalent to asking the fibers $E_b \eqdef \fib p(b)$ to be contractible for all $b: B$, in which case the products in the expression of $F_P$ are indexed over singleton types, hence the name linear.
\end{definition}

\begin{remark}
  \label[remark]{rmk:linear-poly}
  Any span $\begin{tikzcd}[cramped,sep=small]A&\ar[l,"f"']X\ar[r,"g"]&B\end{tikzcd}$ induces a linear polynomial $(X,X,f,\id_X,g)$.
  Moreover this map is an equivalence: every linear polynomial is equivalent (and thus, by univalence, equal) to one of this form.
\end{remark}

\noindent
We can generalize \cref{def:linear-poly} by weakening the requirements on the fibers of $p$, i.e.\ by selecting which \emph{arities} are allowed for the products appearing the associated polynomial functor.

\begin{definition}
    \label[definition]{def:V-ary}
    A polynomial $(E,B,s,p,t)$ is said to be $\cV$\nbd-ary if the fibers of the map $p$ are $\cV$-small.
    We write $\Poly_\cV(I,J)$ for the type of $\cV$-ary polynomials between $I$ and $J$, in other words:
    \[
      \Poly_\cV(I,J) \eqdef \dsum{E: \cU}{\dsum{B: \cU}{(E \to I) \times (E \to_\cV B) \times (B \to J)}}
    \]
    (where $E \to_\cV B$ denotes the type of maps $f : E \to B$ whose fibers are $\cV$-small).
\end{definition}

\begin{example}
    Here are some examples of universes and the corresponding notions of polynomials:
    \begin{itemize}
        \item if $\cV$ is the universe of contractible types, the $\cV$-ary polynomials are the linear ones, hence the spans by \cref{rmk:linear-poly},
        \item if $\cV$ is the universe of propositions (i.e.\ subsingletons), the $\cV$-ary polynomials could be called \emph{affine polynomials},
        \item if $\cV \equiv \Fin$ is the universe of finite types, then we talk about finitary polynomial functors (those correspond to the ones of~\cite{finster2021cartesian}, which are further detailed in \cref{bicategorical-model}).
    \end{itemize}
    Note that those notions make sense even for choices of $\V$ that do not satisfy all the axioms asked at the beginning of \cref{exponential-comonad}.
\end{example}

\begin{proposition}
    \label[proposition]{prop:poly-cat}
    There is a wild category $\Poly_\cV$ whose objects are types in $\cU$ and morphisms are $\cV$\nbd-ary polynomials.
    Identities are given by identity spans seen as polynomials: $\Id_A \eqdef \begin{tikzcd}[cramped,sep=small]A&\ar[l,"\id_A"']A\ar[r,"\id_A"]&A\ar[r,"\id_A"]&A\end{tikzcd}$.
    The composition of two polynomials $P = (\begin{tikzcd}[cramped,sep=small]I&\ar[l,"s"']E\ar[r,"p"]&B\ar[r,"t"]&J\end{tikzcd})$ and $Q = (\begin{tikzcd}[cramped,sep=small]J&\ar[l,"u"']F\ar[r,"q"]&C\ar[r,"v"]&K\end{tikzcd})$ is given by
    \[\begin{tikzcd}
        I & {\sum_{(c,\alpha):D}\sum_{x : F_c} E_{\alpha(x)}} & D & K
        \arrow["{\pi_1}", from=1-2, to=1-3]
        \arrow["{s\circ \pi_3}"', from=1-2, to=1-1]
        \arrow["{v \circ \pi_1}", from=1-3, to=1-4]
      \end{tikzcd}
    \]
    where $D \eqdef \dsum{c : C}{\dprod{x : F_c}B_{u(x)}}$.
\end{proposition}
\begin{proof}
    We do not prove that composition is unital and associative since it will follow from \cref{kleisli-poly}.
    The fact that $\Poly_\V$ is univalent was claimed in \cite{finster2021cartesian} without proof.
    Here is a sketch of proof: suppose that a polynomial $P = (\begin{tikzcd}[cramped,sep=small]I&\ar[l,"s"']E\ar[r,"p"]&B\ar[r,"t"]&J\end{tikzcd})$ has a left and right inverse, we can prove that $p$ has to be an isomorphism in $\V$, so that $P$ is actually a span, and similarly for its two-sided inverses.
    Then univalence of $\Poly_\V$ follows from \cref{span-univalent}.
\end{proof}

\subsection{Polynomials are Kleisli morphisms}
\newcommand{\pts}{\textnormal{poly-to-span}}
\begin{proposition}
  \label[proposition]{prop:poly-span}
  \label[proposition]{poly-span}
    Let $I$, $J$ be types. There is an equivalence
    \[\pts: \Poly_\cV(I,J) \simeq \Span(\bang_\cV I, J)\]
    which maps a polynomial $I \xleftarrow{s} E \xrightarrow{p}_\cV B \xrightarrow{t} J$ to $\bang_\cV I \xleftarrow{\overline s} B \xrightarrow{t} J$ where $\overline s(b) \eqdef \bigl(E_b, \restr{s}{E_b}\bigr)$ (where $\restr{s}{E_b}$ is the restriction of $s$ to~$E_b$).
\end{proposition}
\begin{proof}
    Fix $B: \cU$.
    \begin{equation}
      \label{eq:poly-span}
      \dsum{E: \cU}{(E \to I) \times (E \to_\cV B)}
      \simeq
      \dsum{F: B \to \cV}{\pa{\Sigma B F \to I}}
      \simeq
      B \to \dsum{F: \cV}{(F \to I)}
      \equiv
      (B \to \bang_\cV I)
    \end{equation}
    Thus, we have the following chain of equivalences:
    \begin{align*}
        \Poly_\cV(I,J)
        &\equiv \dsum{E: \cU}{\dsum{B: \cU}{(E \to I) \times (E \to_\cV B) \times (B \to J)}}
            &&\\
        &\simeq \dsum{B: \cU}{\Bigpa{\dsum{E: \cU}{(E \to I) \times (E \to_\cV B)}} \times (B \to J)}
            &&\text{(reordering of terms)}\\
        &\simeq \dsum{B: \cU}{(B \to \bang_\cV I) \times (B \to J)}
            &&\text{(by \cref{eq:poly-span})}\\
        &\equiv \Span(\bang_\cV I, J)
      \end{align*}
      which concludes the proof.
\end{proof}
  
Now, the interesting point is that this equivalence respects composition and identities: as we will show, there is an equivalence of categories between $\Poly$ and the Kleisli category for the comonad $\bang$ on $\Span$. The intuition one should have here is that spans correspond to linear polynomials and $\bang$ allows for using variables many times, we thus expect a linear map which is allowed to use its arguments multiple times to be the same as a polynomial.

\begin{theorem}
  \label[theorem]{kleisli-poly}
  The Kleisli wild category $\Kleisli$ associated to the comonad $\bang_\V$ on $\Span(\U)$ is equivalent to the wild category $\Poly_\V$ of polynomials.
\end{theorem}
\begin{proof}
    The wild categories $\Kleisli$ and $\Poly_\V$ have the same type of objects, namely $\U$, so we can take the identity as the mapping on objects $\Poly_\V \to \Kleisli$.
    Since it is an equivalence on objects, the functor we are constructing is in particular essentially surjective.
    The action on morphisms is given by the equivalence of \cref{poly-span} (thus our functor will be fully faithful).
    Remains to be shown that this mapping is compatible with identities and composition.

    Fix a type $A$.
    By \cref{poly-span}, the identity polynomial on $A$ is mapped to the span $\begin{tikzcd}[cramped,sep=small] \bang_\cV A &\ar[l, "\overline {\id_A}"'] A \ar[r, "\id_A"] & A\end{tikzcd}$, with $\overline{\id_A}(a) = \pa{\dsum{a':A} (a' = a),(a',p) \mapsto a'}$.
    But the type $\dsum{a':A} (a' = a)$ is contractible with center $(a, \refl_a)$~\cite[Lemma~3.11.8]{hottbook}, so $\overline{\id_A}(a) = (\I, \cst_a) = \eta_A(a)$. So in the end $\pts(\Id_A) = \epsilon_A$, which is the identity of $A$ in $\Kleisli$.

    Now for composition. Let $P = (\begin{tikzcd}[cramped,sep=small]I&\ar[l,"s"']E\ar[r,"p"]&B\ar[r,"t"]&J\end{tikzcd}$) and $Q = (\begin{tikzcd}[cramped,sep=small]J&\ar[l,"u"']F\ar[r,"q"]&C\ar[r,"v"]&K\end{tikzcd})$ be polynomials in $\Poly_\V$.
    Writing $D \eqdef \dsum{c : C}{\dprod{x : F_c}B_{u(x)}}$ as in \cref{prop:poly-cat}, the span $\pts(Q \circ P)$ is given by $\begin{tikzcd}[cramped,sep=small] \bang I &\ar[l, "f"'] D \ar[r, "g"] & K\end{tikzcd}$, with
    \[
      f(c,\alpha) \eqdef \Bigl(\dsum{x : F_c} E_{\alpha(x)}, (x,e) \mapsto s(e)\Bigr)
      \qquad\qquad \text{and} \qquad\qquad
      g(c, \alpha) \eqdef v(c)
      \text.
    \]
    On the other hand, composition in $\Kleisli$ is given by the following composition of spans
    \[\begin{tikzcd}[row sep=small]
        & {\bang \bang I} && {\bang B} && C \\
        {\bang I} && {\bang \bang I} && {\bang J} && K
        \arrow["{!\overline s}", from=1-4, to=2-3]
        \arrow["{!t}"', from=1-4, to=2-5]
        \arrow["{\overline u}", from=1-6, to=2-5]
        \arrow["v"', from=1-6, to=2-7]
        \arrow["{\id_{!!I}}"', from=1-2, to=2-3]
        \arrow["\mu", from=1-2, to=2-1]
    \end{tikzcd}\]
    Pulling back along $\id_{\bang \bang I}$ does not change the morphism, so that this is the same as the composition of the spans
    \[\begin{tikzcd}[row sep=small]
        & {\bang B} && C \\
        {\bang I} && {\bang J} && K
        \arrow["{\mu \circ !\overline s}", from=1-2, to=2-1]
        \arrow["{!t}"', from=1-2, to=2-3]
        \arrow["{\overline u}", from=1-4, to=2-3]
        \arrow["v"', from=1-4, to=2-5]
    \end{tikzcd}\]
  This composition is obtained by computing the following pullback:
    \begin{align*}
        {\blue \bang B} {\red \times_{\bang J}} {\green C}
        &\equiv {\blue \dsum{X : \V} \dsum{\pi : X \to B}} {\green\dsum{c : C}} {\red \pa{X , t \circ \pi} =_{\bang J} \pa{F_c , \restr{u}{F_c}}}
            &\\
        &\equivto {\green \dsum{c : C}} {\blue \dsum{X : \V}} {\red \dsum{k : X \equivto F_c}} {\blue \dsum{\pi : X \to B}} {\red t \circ \pi = \restr{u}{F_c} \circ k}
            &&\text{(equality in $\Sigma$-types and univalence)}\\
        &\equivto {\green \dsum{c : C}} {\blue \dsum{\pi : F_c \to B}} {\red t \circ \pi = \restr{u}{F_c}}
            &&\text{(contracting $(X,k)$ onto $(F_c, \id)$)}\\
        &\equivto {\green \dsum{c : C}} {\blue \dsum{\pi : F_c \to B}} {\red \dprod{x : F_c} t(\pi(x)) = u(x)}
            &&\text{(function extensionality)}\\
        &\equivto {\green \dsum{c : C}} {\blue \dprod{x : F_c} \dsum{b : B}} {\red t(b) = u(x)}
            &&\text{(swapping $\Pi$ and $\Sigma$)}\\
        &\equiv {\green \dsum{c : C}} {\color{magenta}\dprod{x : F_c} B_{u(x)}}
            &&\text{(def. of $B_{u(x)} \equiv \fib{t}{u(x)}$)}\\
        &\equiv D &
    \end{align*}
    The inverse map $l : D \to \bang B \times_{\bang J} C$ maps the pair $(c, \alpha)$ to the tuple $((F_c, \fst \circ \alpha),c,(\refl_{F_c}, \snd \circ \alpha))$.
    What remains to show is that the two following triangles commute.
    \[\begin{tikzcd}[row sep=small]
        & D \\
        {\bang I} && J \\
        & {{\blue \bang B} {\red \times_{\bang J}} {\green C}}
        \arrow["l", from=1-2, to=3-2]
        \arrow["f"', from=1-2, to=2-1]
        \arrow["g", from=1-2, to=2-3]
        \arrow["{v \circ {\green \snd}}"', from=3-2, to=2-3]
        \arrow["{\mu \circ \overline s \circ {\blue \fst}}", from=3-2, to=2-1]
    \end{tikzcd}\]
    Let $(c, \alpha): D$. We have
    \begin{align*}
        \mu(\bang \overline s({\blue \fst}(l(\alpha, c)))) &=
        \mu(\bang \overline s(F_c, \alpha)) =
        \mu(F_c, \overline s \circ \alpha) =
        \mu \bigpa{F_c, x \mapsto \bigpa{E_{\alpha(x), \restr{s}{E_{\alpha(x)}}}}}
        = \Bigpa{\dsum{x : F_c} E_{\alpha(x)}, (x,e) \mapsto s(e)}\\
        &=
        f(c, \alpha)
    \end{align*}
    and $v({\green \snd}(l(c, \alpha))) = v(c) = g(c, \alpha)$.
    This completes the proof.
\end{proof}

\noindent
We would like to point out that the observation that polynomial functors could be described as a Kleisli category was already made by Street~\cite[Example~13]{street_polynomials_2020} (although ignoring size issues).
  

\section{Related works}
\label{related-works}

\subsection{Recovering the bicategory of polynomials}
\label{bicategorical-model}
Let us explain here the relationship between the work in this article and the construction of the (wild) bicategory of spans performed in~\cite{finster2021cartesian}.
First, as indicated in \cref{wild-categories}, the notion of wild category is ``right'' only when the hom-types are sets (\cref{category}). In the general case, one needs to add further additional axioms in order to ensure the coherence of the structure (in fact, this coherence is itself part of the structure). The ``next layer'' of coherences can be expressed as follows.


\begin{definition}
  A \emph{$2$-coherence} for a wild precategory~$\C$ consists of witnesses for the commutation of the following diagrams for any suitably composable morphisms $f,g,h,i$:
  \[
    \begin{tikzcd}[row sep=small]
      & {i \circ (h \circ (g \circ f))} \\
      {i \circ ((h \circ g) \circ f)} && {(i \circ h) \circ (g \circ f)} \\
      {(i \circ (h \circ g)) \circ f} && {((i \circ h) \circ g) \circ f}
      \arrow["{\alpha_{i \circ h, g, f}}", from=2-3, to=3-3]
      \arrow[""{name=0, anchor=center, inner sep=0}, "{\ap{(- \circ f)}(\alpha_{i,h,g})}"', from=3-1, to=3-3]
      \arrow["{\alpha_{i,h\circ g,f}}"', from=2-1, to=3-1]
      \arrow["{\ap{(i\circ-)}(\alpha_{h,g,f})}"', from=1-2, to=2-1]
      \arrow["{\alpha_{i,h,g\circ f}}", from=1-2, to=2-3]
    \end{tikzcd}
    \hfill
    \begin{tikzcd}[row sep=small]
      {g \circ (\id \circ f)} && {(g \circ \id) \circ f} \\
      & {g \circ f}
      \arrow[""{name=0, anchor=center, inner sep=0}, "{\alpha_{g,\id,f}}", from=1-1, to=1-3]
      \arrow["{\ap{(g \circ -)}(\lambda_f)}"', from=1-1, to=2-2]
      \arrow["{\ap{(- \circ f)}(\rho_g)}", from=1-3, to=2-2]
    \end{tikzcd}
  \]
  %
  A (\emph{pre})\emph{bicategory} is a wild (pre)category equipped with a $2$-coherence such that $\C(x,y)$ is a groupoid for every objects $x,y:\C$.
\end{definition}

It can be shown that the wild category $\Span(\U)$ can be equipped with a $2$-coherence. We expect that the monad defined in \cref{exponential-comonad} can also be shown to be coherent (this is left for future work), thus inducing a $2$-coherence on the wild Kleisli category $\Poly(\U)$. Moreover, this construction can be refined (in the same spirit as \cref{small-type}) in order to replace $\U$ by a subuniverse of types satisfying suitable closure properties. In particular, this would allow defining $\Poly(\Gpd)$, thus recovering the construction of the wild bicategory of polynomials defined in~\cite{finster2021cartesian}. It should be noted that, contrarily to what is claimed, \cite{finster2021cartesian} only actually constructs a wild $2$-coherent category: it is $3$-truncated (and not $2$-truncated as expected for a bicategory). However, this can be fixed by truncating the wild $2$-coherent category into a bicategory, using a similar process as for categories (see \cref{wild-category-truncation,wild-category-truncation-univalence}).

\subsection{Towards a model of differential linear logic}
\label{template-games}
In~\cite{mellies_template_2019}, Melliès defines a bicategorical model of differential linear logic based on a comonad on spans (we consider here only the particular case where the ``synchronization template'' category is the terminal one). When we restrict his model to groupoids (as opposed to categories), we obtain a model which is equivalent to the one of the previous section with $\U\defd\Gpd$ and $\V\defd\Fin$. We can thus conclude that the morphisms in his Kleisli category are the polynomials, this observation being new to our knowledge. This also suggests that we could extend our model into one of differential linear logic. This shall be the topic of future work.

\subsection{Related models}
The model we have defined is close to the one of species of structures~\cite{fiore2008cartesian} which are a categorical generalization of Joyal's combinatorial species~\cite{joyal1981theorie}. We plan to investigate the relationship between the two in future work. In particular, the comparison should be made easier by the work of Gepner, Haugseng and Kock~\cite[Section 3.2]{gepner2022operads} who showed that analytic and polynomial functors coincide in the setting of $\infty$-categories.
It would also be interesting to compare the present work with Clairambault and Forest model of thin spans of groupoids~\cite{clairambault2023cartesian}, which is close in spirit but technically quite different since it uses an orthogonality construction in order to avoid taking symmetries in account.

Some other authors have considered categories whose objects (as opposed to morphisms) are polynomials. Such categories have also been shown to be cartesian closed~\cite{altenkirch2010higher}, to define models of linear logic~\cite{hyvernat}, and to be the Kleisli category of some other category~\cite[Theorem~2.4.5]{kocktrees}. Although those results seem reminiscent of the ones exposed here, it is unclear if they can be related in any way, since there is in general no connection between a bicategory and the category of endomorphisms of one of its objects.

\bibliographystyle{entics}
\bibliography{biblio,samuel}


\end{document}